\documentclass{ansarticle}

\usepackage[utf8]{inputenc}
\usepackage{subcaption}

\newcommand{\dune}{\textsc{Dune}}

\newcommand{\dunetypetree}{\textsc{Dune-TypeTree}}
\newcommand{\dunefunctions}{\textsc{Dune-Functions}}
\newcommand{\dunedpg}{\textsc{Dune-DPG}}

\newcommand{\cB}{\ensuremath{\mathcal{B}}}
\newcommand{\cC}{\ensuremath{\mathcal{C}}}

\newcommand{\cN}{\ensuremath{\mathcal{N}}}

\newcommand{\bN}{\ensuremath{\mathbb{N}}}

\newcommand{\bP}{\ensuremath{\mathbb{P}}}

\newcommand{\bR}{\ensuremath{\mathbb{R}}}

\newcommand{\bU}{\ensuremath{\mathbb{U}}}
\newcommand{\bV}{\ensuremath{\mathbb{V}}}

\DeclareMathOperator*{\vspan}{span}

\newcommand{\supp}{\ensuremath{\mathrm{supp}}}

\newcommand{\opt}{\ensuremath{\mathrm{opt}}}
\newcommand{\no}{\ensuremath{\mathrm{n.opt}}}

\newcommand\restr[2]{{
  \left.\kern-\nulldelimiterspace 
  #1 
  \vphantom{\big|} 
  \right|_{#2} 
}}

\title{
The \dunedpg~library for solving PDEs with Discontinuous Petrov--Galerkin finite elements
}
\author[1]{Felix Gruber}
\author[1]{Angela Klewinghaus}
\author[2]{Olga Mula}
\affil[1]{IGPM, RWTH Aachen, Templergraben 55, 52056 Aachen, Germany}
\affil[2]{Université Paris-Dauphine, PSL Research University, CEREMADE, 75775 Paris, France}
\begin{document}

\maketitle

\begin{abstract}
In the numerical solution of partial differential equations (PDEs), a central question is the one of building variational formulations that are inf-sup stable not only at the infinite-dimensional level, but also at the finite-dimensional one. This guarantees that residuals can be used to tightly bound errors from below and above and is crucial for a posteriori error control and the development of adaptive strategies. In this framework, the so-called Discontinuous Petrov--Galerkin (DPG) concept can be viewed as a systematic strategy of contriving variational formulations which possess these desirable stability properties, see e.\,g.\ \citet{BDS2015}.
In this paper, we present a C++ library, \dunedpg, which serves to implement and solve such variational formulations. The library is built upon the multi-purpose finite element package \dune~(see \citet{dune2.4}). One of the main features of \dunedpg\ is its flexibility which is achieved by a highly modular structure. 
The library can solve in practice some important classes of PDEs (whose range goes beyond classical second order elliptic problems and includes e.\,g.\ transport dominated problems). As a result, \dunedpg~can also be used to address other problems like optimal control with the DPG approach.
\end{abstract}

\section{Introduction}
\paragraph*{General context and motivations:} Let $\Omega$ be a domain of $\bR^d$ ($d\geq 1$) and $\bU,\ \bV$ two Hilbert spaces defined over $\Omega$ and endowed with norms $\Vert \cdot \Vert_{\bU}$ and $\Vert \cdot \Vert_{\bV}$, respectively. The normed dual of $\bV$, denoted $\bV'$, is endowed with the norm
\begin{equation*}
\Vert \ell \Vert_{\bV'} \coloneqq \sup_{v\in \bV} \frac{\vert \ell(v) \vert}{\Vert v \Vert_{\bV}},\quad \forall \ell \in \bV'.
\end{equation*}
Let $\cB:\bU\to \bV$ be a boundedly invertible linear operator and let $b:\bU\times \bV \to \bR$ be its associated continuous bilinear form defined by $b(w,v)=\left( \cB w\right)(v),\ \forall (w,v)\in \bU\times \bV$. We consider the operator equation
\begin{align}
\label{eq:operatorEq}
\begin{split}
&\text{Given $f\in \bV'$, find $u\in \bU$ s.\,t.} \\
&\cB u = f,
\end{split}
\end{align}
or, equivalently, the variational problem
\begin{align}
\label{eq:WFinfite}
\begin{split}
&\text{Given $f\in \bV'$, find $u\in \bU$ s.\,t.} \\
&b(u,v) = f(v), \quad \forall v \in \bV.
\end{split}
\end{align}
Let $0<\gamma \leq 1$ be a lower bound for the (infinite-dimensional) inf-sup constant
\begin{equation*}
\inf_{w\in \bU} \sup_{v\in \bV} \frac{b(w,v)}{\Vert w \Vert_{\bU} \Vert v\Vert_{\bV}} \geq \gamma >0.
\end{equation*}
Since $\cB$ is invertible, problem \eqref{eq:operatorEq} admits a unique solution $u\in \bU$ and for any approximation $\bar u \in \bU$ of $u$,
\begin{equation}
\label{eq:resBoundInfinite}
\Vert \cB \Vert_{L(\bU,\bV')}^{-1} \Vert f - \cB \bar u \Vert_{\bV'}
\leq
\Vert u - \bar u \Vert_{\bU}
\leq
\gamma^{-1} \Vert f - \cB \bar u \Vert_{\bV'}.
\end{equation}
From \eqref{eq:resBoundInfinite}, it follows that the error $\Vert u - \bar u \Vert_{\bU}$ is equivalent to the residual $\Vert f - \cB \bar u \Vert_{\bV'}$. The residual contains known quantities and its estimation on an appropriate finite-dimensional space opens the door to rigorously founded a posteriori concepts.
However, note that the information that the estimator can give is only meaningful when the variational formulation is well-conditioned, i.\,e., for $\Vert \cB \Vert_{L(\bU,\bV')}$ and $\gamma$ being as close to one as possible.
Assuming that we have this property of well-conditioning, a crucial point is that this needs to be inherited at the finite-dimensional level.
This issue has been well explored for standard Galekin methods (i.\,e.~when $\bU=\bV$) and allows to appropriately address most parabolic and second order elliptic problems with a wide variety of finite element methods.
However, much less is known when $\bU\neq\bV$ is required to obtain a well-posed and well-conditioned problem, like for transport-dominated PDEs. In the latter case, the main ideas are that
\begin{itemize}
\item $\gamma=1$ by choosing a problem-dependent norm for $\bV$ and 
\item for a given finite-dimensional trial space $\bU_H$, there exists a corresponding optimal test space $\bV^{opt}(\bU_H)$ such that the discrete inf-sup condition
\begin{equation}
\inf_{w_{H}\in \bU_H}
\sup_{v\in \bV^{\opt}(\bU_H)}
\frac{b(w_H,v)
}{
\Vert w_{H} \Vert_{\bU} \Vert v\Vert_{\bV}
}
\geq \gamma
\end{equation}
holds with the same constant $\gamma$ as the infinite-dimensional one. For more details on this, see Section~\ref{sec:optTestSpacesTheory}.
\end{itemize}
Since, in general, even the approximate computation of the optimal test space requires the solution of global problems, there are essentially two ways to make the computation affordable. One is the introduction of a mixed formulation  which avoids the computation of the optimal test spaces and only uses them indirectly. This approach was used in \citet{DHSW2012} to construct a general adaptive scheme when $\bU = L_2(\Omega)$ and to show convergence under certain abstract conditions (which have to be verified for concrete applications). It was also employed in the context of reduced-basis construction for transport-dominated problems (see \citep{DPW2014}).

The other way is to make computations affordable by localization so that the optimal test spaces can be computed by solving local problems. This is the approach taken in the DPG methodology, initiated and developed mainly by L.~Demkowicz and J.~Gopalakrishnan (see e.\,g.~\citet{DG2011a,GQ2014,DG2015}). In this method, an approximation of the exact optimal test functions is realized in the context of a discontinuous Petrov--Galerkin formulation. This yields a so-called near-optimal test space which is the one that is eventually used in the solution of the discrete problem.
The tightness of a posteriori error estimators has been theoretically justified only for second order elliptic problems.
However, numerical evidence illustrates their good performance also for a much broader variety of problems.
Without being exhaustive, we can find works on transport equations~\citep{BDS2015}, convection--diffusion~\citep{BS2015}, elasticity and Stokes problems~\citep{CDG2014}, Maxwell equations~\citep{CDG2015} and the Helmholtz equation \citep{DGMZ2012}.

The presented finite element library \dunedpg\ serves to implement and solve such DPG formulations.

\paragraph*{Contributions and layout of the paper:}
In this paper, we explain the construction of the \dunedpg~library which is capable of solving various types of PDEs with a DPG variational formulation. Since the appropriate characteristics of the formulation depend on the problem, the user is given the freedom to choose the spaces $\bU$, $\bV$, their finite-dimensional counterparts and the geometry and mesh refinement.
We would like to note that, in fact, \dunedpg\ is not the first software tool for solving PDEs with the DPG method.
The \emph{Camellia} package (see \citet{Roberts2014}) is another library whose purpose and construction is similar to \dunedpg.
In \dunedpg, we incorporate, at the software level, the latest theoretical results on DPG for transport equations, given by \citet{BDS2015}, that allow to appropriately address families of transport-based PDEs.
They essentially require the use of subgrids of fixed depth for the finite-dimensional test space.
For this reason, in \dunedpg, a strong emphasis has been put on providing as much freedom as possible in the selection of the finite-dimensional test space.
This makes \dune\ an appropriate choice of foundation for our library since its modular structure gives low-level access to those parts for which fine control is required while still providing high-level functionality for the rest of the code.

To show how the library works, the paper is organized as follows: in Section~\ref{sec:foundationsDPG}, we summarize the mathematical concepts of DPG that are relevant to understand the library. As an example, we explain at the end of this section how the ideas can be applied to a simple transport problem following the theory of \citet{BDS2015}. Then, in Section~\ref{sec:dunedpg}, we present the different building blocks that form the library. We explain how they interact and how they make use of some features of the \dune\ framework upon which our library is built. Finally, in Section~\ref{sec:numEx}, we validate \dunedpg\ by giving concrete results related to the solution of a simple transport problem.

\section{Theoretical Foundations for DPG}
\label{sec:foundationsDPG}
As already brought up in the introduction, the DPG concept was initiated and developed mainly by L.~Demkowicz and J.~Gopalakrishnan (see e.\,g.~\citet{DG2011a,GQ2014}). Other relevant results concerning theoretical foundations are \citet{BS2014b,BS2015} and, more recently, \citet{BDS2015}. The strategy followed in DPG to contrive stable variational formulations is based on the concept of \textit{optimal test spaces} and their practical approximation through the solution of \textit{local} problems in the context of a discontinuous Petrov--Galerkin variational formulation. The two following sections explain more in detail these two fundamental ideas.

\subsection{The concepts of optimal and near-optimal test spaces}
\label{sec:optTestSpacesTheory}
Assuming that we start from a well-posed and well-conditioned infinite-dimensional variational formulation \eqref{eq:WFinfite}, we look for a formulation at the finite-dimensional level which inherits these desirable features. Let $H>0$ be a parameter ($H$ will later be associated to the size of a mesh $\Omega_H$ of $\Omega$). For any given finite-dimensional trial space $\bU_{H}$ of dimension $\cN$ (that depends on $H$), there exists a so-called optimal test space $\bV^{\opt}(\bU_{H})$ of the same dimension. It is called optimal because the finite-dimensional version of problem \eqref{eq:WFinfite},
\begin{equation}
\begin{aligned}
&\text{Find $u_{H}\in \bU_{H}$ s.\,t.} \\
&b(u_{H},v) = f(v), \quad \forall v \in \bV^{\opt}(\bU_{H}),
\end{aligned}
\label{eq:WFfiteOpt}
\end{equation}
is well posed and
\begin{equation*}
\inf_{w_{H}\in \bU_H}
\sup_{v\in \bV^{\opt}(\bU_H)}
\frac{b(w_H,v)
}{
\Vert w_{H} \Vert_{\bU} \Vert v\Vert_{\bV}
}
\geq \gamma .
\end{equation*}
In other words, the discrete inf-sup condition is bounded with the same constant $\gamma$ that is involved in the infinite-dimensional problem. This implies that the discrete problem has the same stability properties as the infinite-dimensional problem. Therefore the residual $\Vert f - \cB u_{H} \Vert_{\bV'}$ is equivalent to the actual error $\Vert u - u_{H}\Vert_{\bU}$ with the same constants exhibited in \eqref{eq:resBoundInfinite}. Since these constants do not depend on $H$, $\Vert f - \cB u_{H} \Vert_{\bV'}$ is a \textit{robust} error bound that is suitable for adaptivity since we can decrease $H$ without degrading the constants of equivalence.

Unfortunately, the optimal test space $\bV^{\opt}(\bU_H)$ is not computable in practice. Indeed, if $\{u_H^i\}_{i=1}^{\cN}$ spans a basis of $\bU_H$, then the set of functions $\{v^i\}_{i=1}^{\cN}$ defined through the variational problems,
\begin{align}
\label{eq:optBasisInfinite}
i\in \{1,\dots,\cN\},\quad
\langle v^i,v \rangle_{\bV} = b(u_H^i,v),\quad \forall v \in \bV
\end{align}
spans a basis of $\bV^{\opt}(\bU_H)$. Since these problems are formulated in the infinite-dimensional space $\bV$, they cannot be computed exactly (in addition, the problems are global). 
To address this issue, problems \eqref{eq:optBasisInfinite} are $\bV$-projected to a finite-dimensional subspace $\bV_h$ that will be called test-search space.
Therefore, in practice, an approximation $\{\bar v^i\}_{i=1}^{\cN}$ to the set of functions $\{v^i\}_{i=1}^{\cN}$ is computed by solving for all $i\in \{1,\dots,\cN\}$,
\begin{align}
\label{eq:optBasisFinite}
\langle \bar v^i,v \rangle_{\bV} = b(u_H^i,v),\quad \forall v \in \bV_h.
\end{align}
This defines a projected test space $\bV^{\no}(\bU_H,\bV_h)\coloneqq\vspan\{\bar v^i\}_{i=1}^{\cN}$. For the elliptic case and some classes of transport problems, it is possible to exhibit test-search spaces $\bV_h$ (which depend on the initial $\bU_H$) such that $\bV^{\no}(\bU_H,\bV_h)$ is close enough to the optimal $\bV^{\opt}(\bU_H)$ to allow that the discrete $\inf$-$\sup$ constant
\begin{equation}
\label{eq:condAdaptivity}
\gamma_{H} \coloneqq
\inf_{u_H\in \bU_{H}} \sup_{v_h\in \bV^{\no}(\bU_H,\bV_h)} \frac{b(u_H,v_h)}{\Vert u_H \Vert_{\bU_H} \Vert v_h\Vert_{\bV_h}}
\end{equation}
is bounded away from $0$ uniformly in $H$. For this reason, $\bV^{\no}(\bU_H,\bV_h)$ is called a near-optimal test-space. In the case of transport problems, the recent work of \citet{BDS2015} shows that good test-search spaces $\bV_h$ can be found when they are defined over a refinement $\Omega_h$ of $\Omega_H$.

The near-optimal test space $\bV^{\no}(\bU_H,\bV_h)$ is the one that is computed in practice in the \dunedpg~library. The finite-dimensional variational formulation that is eventually solved reads
\begin{align}
\label{eq:WFfiteNearOpt}
\begin{split}
&\text{Find $u_H \in \bU_H$ s.\,t.} \\
&b(u_H,v_h) = f(v_h), \quad \forall v_h \in \bV^{\no}(\bU_H,\bV_h).
\end{split}
\end{align}
It can be expressed as a linear system of the form $Ax=F$, $A \in \bR^{\cN\times \cN},\ x\in \bR^{\cN},\ F\in \bR^{\cN}$. It can be proven that $A$ is by construction symmetric positive definite. The assembly of the system and its solution in \dunedpg\ are explained in Section~\ref{sec:assembly}.

\subsection{The concept of localization}
\label{sec:localizationTheory}
Depending on the choice of $\bV$ and $\bV_h$, the solution of \eqref{eq:optBasisFinite} to derive the near-optimal basis functions of $\bV^{\no}(\bU_H,\bV_h)$ might be costly. This is because these $\cN$ problems are, in general, global in the whole domain $\Omega$ and they cannot be decomposed into local ones. Furthermore, if the resulting near-optimal basis functions have global support, the solution of the finite-dimensional variational problem \eqref{eq:WFfiteNearOpt} is costly as well because the resulting system matrix $A$ is full.

To prevent this, we need an appropriate variational formulation with a well-chosen test space $\bV$ which has a product structure on the coarse grid $\Omega_H$,
\begin{equation}
\label{eq:Vprod}
\bV \coloneqq \prod_{K \in \Omega_H} \bV_K,
\end{equation}
where $\supp(v) \subset K$ for any $v \in \bV_K$. In particular, the restriction of the $\bV$-scalar product to $K \in \Omega_H$ has to be a scalar product for $\bV_K$:
\begin{equation*}
\langle \cdot, \cdot \rangle_{\bV} \vert_K = \langle \cdot, \cdot \rangle_{\bV_K}
\end{equation*}
The test-search space $\bV_h$ will be choosen in such a way, that it has the same product structure as $\bV$,
\begin{equation*}
\bV_h \coloneqq \prod_{K \in \Omega_H} \bV_{h,K},\quad \bV_{h,K} \subset \bV_K.
\end{equation*}
Therefore, for any $1\leq i \leq \cN$, the near-optimal test function $\bar v^i$ can be written as
\begin{equation*}
\bar v^i = \sum_{K \in \Omega_H} \bar v^i_K \chi_K,
\end{equation*}
where $\chi_K$ is the characteristic function of cell $K$.
Additionally, we need a decomposition of the bilinear form as a sum over mesh cells of $\Omega_H$,
\begin{equation}
\label{eq:bilinSumK}
b\left( u,v \right) = \sum_{K\in \Omega_H} b_{K}\left( u,v \right),\quad \forall v \in \prod_{K \in \Omega_H} \bV_K.
\end{equation}
Then, for every $K\in \Omega_H$, $\bar v^i_K$ is the solution of a local problem in $K$,
\begin{equation}
\label{eq:optBasisFiniteLocal}
\langle \bar v^i_K,v \rangle_{\bV_K} = b_{K} \left( u_H^i,v \right),\ \forall v \in \bV_{h,K},
\end{equation}
where $\{u^i_H\}_{i=1}^{\cN}$ is a basis of $\bU_H$. Therefore, finding $\bar v^i$ can be decomposed into a sum of problems, each one of which is localized on a mesh cell $K\in \Omega_H$. Moreover, if the support of $u^i_H$ is included in some cell $K \in \Omega_H$, then the support of its corresponding near-optimal test function $\bar v^i$ is also a subset of $K$ (and the neighboring cells in some cases). In other words, we would have $\bar v^i = \bar v^i_K \chi_K$ or $\bar v^i = \sum_{K' \in\ \text{neigh(K)}} \bar v^i_{K'} \chi_{K'}$. Hence, if the basis functions $u^i_H$ of $\bU_H$ have local support, the resulting system matrix $A$ is sparse.

\subsection{An example: a linear transport equation}
\label{sec:introTransport}
Let $\Omega = (0,1)^2$ and $\beta$ be a vector of $\bR^2$ with norm one. For any $x\in \partial \Omega$, let $n(x)$ be its associated outer normal vector. Then
\begin{equation}
\label{Gammas}
\Gamma_- \coloneqq \{ x \in \partial \Omega  \mid  \beta\cdot n(x) < 0 \} \subset \partial\Omega
\end{equation}
is the inflow-boundary for the given constant transport direction $\beta$. Given $c\in \bR$ and a function $f:\Omega \to \bR$, we consider the problem of finding the solution $\varphi:\Omega \to \bR$ to the simple transport equation
\begin{equation}
\label{eq:pureTransport}
\begin{alignedat}{2}
\beta \cdot \nabla \varphi + c \varphi &= f,&& \quad \text{in $\Omega$},\\
\varphi &= 0,&& \quad \text{on $\Gamma_-$.}
\end{alignedat}
\end{equation}
If we apply the DPG approach introduced in~\citet{BDS2015} to solve this problem, we first need to introduce the following spaces. Denoting $\nabla_H$ the piecewise gradient operator, let
\begin{equation*}
H(\beta,\Omega_H) \coloneqq \{ v\in L_2(\Omega) \mid \beta\cdot\nabla_H v \in L_2(\Omega)\},
\end{equation*}
equipped with squared ``broken'' norm $\Vert v \Vert^2_{H(\beta,\Omega_H)} = \Vert v \Vert^2_{L_2(\Omega)}+ \Vert \nabla_H\cdot v \Vert^2_{L_2(\Omega)}$. Let also
\begin{equation*}
H_{0,\Gamma_-}(\beta,\Omega) \coloneqq \text{clos}_{H(\beta,\Omega)} \{ u\in H(\beta,\Omega) \cap \cC(\bar \Omega) \mid u=0 \text{ on }\Gamma_- \}
\end{equation*}
and
\begin{equation*}
H_{0,\Gamma_-}(\beta,\partial \Omega_H) \coloneqq \{ w\vert_{\partial \Omega_H} \mid w \in H_{0,\Gamma_-}(\beta,\Omega) \}
\end{equation*}
equipped with quotient norm
\begin{equation*}
\Vert \theta \Vert_{H_{0,\Gamma_-}(\beta,\partial \Omega_H)}
\coloneqq
\inf \{
\Vert w \Vert_{H(\beta,\Omega)} \mid \theta = w\vert_{\partial \Omega_H},\ w\in H_{0,\Gamma_-}(\beta,\Omega)
\}.
\end{equation*}
The variational formulation reads
\begin{equation}
\label{eq:infiniteVFtrans}
\begin{aligned}
&\text{For $\bU \coloneqq L^2(\Omega)\times H_{0,\Gamma_-}(\beta,\partial \Omega_H)$ and $\bV\coloneqq H(\beta,\Omega_H)$,} \\
&\text{given $f \in H(\beta,\Omega_H)'$, find $u \coloneqq (\varphi,\theta)\in \bU$ such that} \\
&b(u,v)=f(v),\quad \forall v \in \bV.
\end{aligned}
\end{equation}
In this formulation (usually called \textit{ultra-weak} formulation) the bilinear form $b(u,v)$ is defined by
\begin{equation}
b(u,v)= b\left( (\varphi,\theta),v\right)
=
\int_{\Omega} \left( -\beta\cdot \nabla v \varphi + cv\varphi \right)\,\dx + \int_{\partial \Omega_H} \llbracket v\beta \rrbracket \theta\,\ds.
\label{eq:ultraweak}
\end{equation}
Note that this variational formulation depends on the mesh $\Omega_H$. Also, note the presence of an additional unknown $\theta$ that lives on the skeleton $\partial \Omega_H$ of the mesh. For smooth solutions, $\theta$ agrees with the traces of $\varphi$ on $\partial \Omega_H$ (i.\,e.~the union of cell interfaces of $\Omega_H$).

For the discretization, we take for some $m\in\bN$,
\begin{equation}
\label{Uh}
\bU_H\coloneqq \bigg( \prod_{K\in\Omega_H}\bP_{m-1}(K)\bigg) \times
\restr{\bigg(H_{0,\Gamma_-}({ \beta};\Omega) \cap
       \prod_{K\in\Omega_H}\bP_m(K)\bigg)}{\partial \Omega_h},
\end{equation}
where $\bP_m(K)$ is the space of polynomials of degree $m$. A viable test-search space can be taken simply as discontinuous piecewise polynomials of slightly higher degree on the finer mesh $\Omega_h$ of $\Omega_H$, namely
\begin{equation}
\label{testsearch}
\bV_h \coloneqq \prod_{K\in \Omega_h}\bP_{m+1}(K).
\end{equation}
As a result, the discrete version of \eqref{eq:infiniteVFtrans} reads
\begin{equation}
\label{eq:finiteVFtrans}
\begin{aligned}
&\text{Find $u_H \coloneqq (\varphi_H,\theta_H)\in \bU_H$ such that} \\
&\tilde b(u_H,v_h)=f(v_h),\quad \forall v_h \in \bV^{\no}(\bU_H,\bV_h).
\end{aligned}
\end{equation}
The bilinear form $\tilde b$ is slightly different from $b$. It reads
\begin{equation}
\tilde b(u,v)= b\left( (\varphi,\theta),v\right)
=
\int_{\Omega} \left( -\beta\cdot \nabla v \varphi + cv\varphi \right)\,\dx + \int_{\partial \Omega_h} \llbracket v\beta \rrbracket \theta\,\ds,
\label{eq:ultraweakdiscrete}
\end{equation}
where the trace integral is over $\partial\Omega_h$ and not $\partial\Omega_H$.

\section{An Overview of the Architecture of \dunedpg}
\label{sec:dunedpg}
In this section we describe how the DPG method presented in Section~\ref{sec:foundationsDPG} has been implemented in \dunedpg. As already brought up, the present library has been built upon the finite element package \dune. It benefits from the new \dunefunctions\ module~\citep{EGMS2015} that has been critical for the construction of a uniform interface for test and trial spaces.

\begin{figure}[h]
\begin{center}
\begin{tikzpicture}[sibling distance=15em,
  every node/.style = {shape=rectangle, rounded corners,
    draw=black, align=center,inner sep=6pt,
    top color=white, bottom color=gray!10}]
  \node (sa) {\lstinline[style=cppstyle]{SystemAssembler}};
  \node (ot) [below left = 10mm and -3mm of sa] {\lstinline[style=cppstyle]{OptimalTestBasis} \\ \footnotesize{for $\bV^{\no}(\bU_H,\bV_h)$}};
  \node (rhs)[right = 15mm of ot] {right hand side\\ \footnotesize{implemented as}\\ \footnotesize{tuple of functions}};
  \node (bc) [right = 3mm of rhs] {boundary conditions\\ \footnotesize{implemented as}\\ \lstinline[style=cppstyle]{applyDirichletBoundarySolution()}};
  \node (ip) [below left = 10mm and -13mm of ot] {\lstinline[style=cppstyle]{InnerProduct}\\ \footnotesize{$\langle \cdot, \cdot \rangle_{\bV_K}$}}; 
  \node (bf) [right = 10mm of ip] {\lstinline[style=cppstyle]{BilinearForm}\\ \footnotesize{$b_K(\cdot, \cdot)$}};
  \node (it) [below left = 10mm and -13mm of ip] {\lstinline[style=cppstyle]{IntegralTerm}s\\ \footnotesize{of $\langle \cdot, \cdot \rangle_{\bV_K}$}};
  \node (vh) [right = 3mm of it] {enriched test space\\ \footnotesize{$\bV_h$}};
  \node (uh) [right = 3mm of vh] {trial space\\ \footnotesize{$\bU_H$}};
  \node (itb)[right = 3mm of uh] {\lstinline[style=cppstyle]{IntegralTerm}s\\\footnotesize{ of $b_K(\cdot, \cdot)$}};

  \draw (sa) -- (ot);
  \draw (sa) -- (bf);
  \draw (sa) -- (rhs);
  \draw (sa) -- (bc);
  \draw (ot) -- (ip);
  \draw (ot) -- (bf);
  \draw (ip) -- (it);
  \draw (ip) -- (vh);
  \draw (bf) -- (vh);
  \draw (bf) -- (uh);
  \draw (bf) -- (itb);
\end{tikzpicture}
\end{center}
\caption{Overview of interactions between the main classes of \dunedpg.}
\label{fig:dunedpgOverview}
\end{figure}
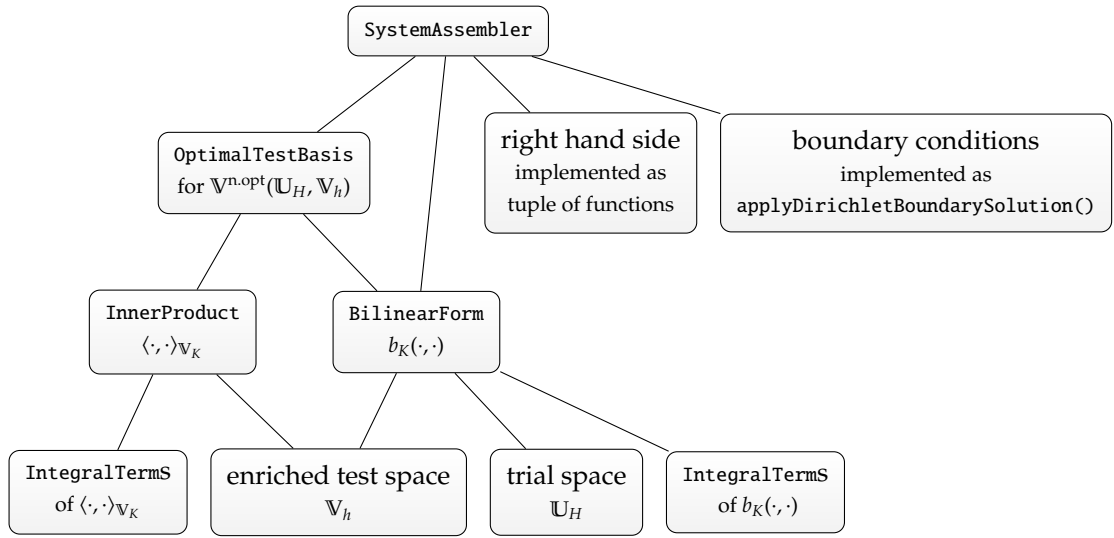

The user of \dunedpg\ starts by choosing the appropriate test-search space $\bV_h$ and trial space $\bU_H$ for his problem. Then, the bilinear form $b(\cdot, \cdot)$ and the inner product $\langle \cdot, \cdot \rangle_{\bV}$ are declared via the classes \cpp{BilinearForm} and \cpp{InnerProduct} (see Section~\ref{sec:bilinform}). They both consist of an arbitrary number of elements of the type \cpp{IntegralTerm} (see Section~\ref{sec:integralterm}). Next, the near-optimal test space is determined automatically with the help of the bilinear form and the inner product via \eqref{eq:optBasisFiniteLocal} (see Section~\ref{sec:optTestSpacesImpl}). Finally, the \cpp{SystemAssembler} class handles the automatic assembly of the linear system $Ax=F$ associated to problem \eqref{eq:WFfiteNearOpt} including the right hand side $f$ and boundary conditions (see Section~\ref{sec:systemassembler}). For assembling, the matrix $A$, we use \eqref{eq:bilinSumK} and define the local matrices $A_K$ by
\begin{equation}
\label{eq:AK}
(A_K)_{i,j}=b_K(u_{H,K}^i, \tilde{v}_K^j)\,
\end{equation}
where $\{ u_{H,K}^i \}_{i=1}^{\mathcal{N}_K}$ is a basis for the restriction of $\bU_H$ to $K$ and $\{\tilde{v}_K^j\}_{j=1}^{\mathcal{N}_K}$ is a basis for the restriction of the near-optimal test spaces $\bV^{\no}(\bU_H,\bV_h)$ to $K$. The \cpp{BilinearForm} class provides the local matrices $A_K$ and the  \cpp{SystemAssembler} class constructs the global matrix $A$ out of the local matrices $A_K$. The classes that we have just mentioned are intertwined and depend on each other. Figure \ref{fig:dunedpgOverview} gives an overview of their interactions. In addition to these classes, the class \cpp{ErrorTools} handles the computation of a posteriori estimators following the guidelines that are given in Section~\ref{sec:apost}.

\subsection{Assembling the discrete system for a given PDE}
\label{sec:assembly}
The following subsections describe the \cpp{SystemAssembler} class and all the classes used by it, except for the \cpp{OptimalTestBasis} class that will be explained in detail in Section~\ref{sec:optTestSpacesImpl}.

\subsubsection{SystemAssembler}
\label{sec:systemassembler}

The assembly of the discrete system $Ax=F$ derived from the variational problem is handled by the class \cpp{SystemAssembler}. We start by creating the appropriate object of the class (that we will name in our case \cpp{systemAssembler}) by calling the method \cpp{make_DPG_SystemAssembler}. As an input, it needs objects representing our trial space $\bU_H$, our near-optimal test space $\bV^{\no}(\bU_H,\bV_h)$ and the bilinear form $b$. The bilinear form is an object of the class \cpp{BilinearForm} which is explained in Section~\ref{sec:bilinform}. As for the spaces $\bU_H$ and $\bV^{\no}(\bU_H,\bV_h)$, they are respectively given by a \cpp{std::tuple} composed of global basis functions from \dunefunctions. The reason to use a tuple is to handle problems involving several unknowns. For instance, in the ultra-weak formulation introduced for the transport problem in Section~\ref{sec:introTransport}, we have two unknowns $(\varphi,\theta)$ and $\bU_H$ is a product of two spaces.

Once \cpp{systemAssembler} has been defined, a call to the method \cpp{assembleSystem(stiffnessMatrix, rhsVector, rhsFunction)} assembles the matrix $A$ and the right-hand side vector $F$. They are stored in the variables \cpp{stiffnessMatrix} (of type \cpp{BCRSMatrix<FieldMatrix<double,1,1> >}) and \cpp{rhsVector} (of type \cpp{BlockVector<FieldVector<double,1> >}). The input parameter \cpp{rhsFunction} is a \cpp{std::tuple} of \cpp{std::function<double(Dune::FieldVector<double,dim>)>} and represents the function $f$ from the PDE. Internally, the class \cpp{SystemAssembler} iterates over all mesh cells $K$ and delegates the work of computing local contributions $A_K$ to the system matrix $A$ to \cpp{BilinearForm::getLocalMatrix()}. Similarily the local right-hand side vectors are computed by a function \cpp{getVolumeTerm}. For constructing $A$ and $F$ out of the local matrices $A_K$ and the local right-hand side vectors, we make use of the mapping between local and global degrees of freedom given by the index sets from \dunefunctions.

Once $A$ and $F$ are obtained, the system $Ax=F$ (which is, from the theory, invertible) can be solved with the user's favorite direct or iterative scheme. As a summary, the following lines of code give the main guidelines to solve the transport problem of Section~\ref{sec:introTransport}. The commented lines starting with ``[...]'' mean that there is some code to be written in addition. We do not provide it for the sake of clarity and refer to our source code for an example of exact implementation (see \lstinline[style=anycodestyle,basicstyle=\ttfamily]{src/plot_solution.cc}).

\lstset{escapeinside={(*}{*)},style=cppstyle}
\begin{c++}
// Definition of the trial spaces
DuneFEM1 spacePhi;
DuneFEM2 spaceTheta;
auto solutionSpaces = std::make_tuple(spacePhi, spaceTheta);
// [...] Definition of testSearchSpaces (tuple containing the test search space)
// [...] Computation the associated near optimal test space. The space is stored in the object nearOptTestSpaces (see Section(*~\ref{sec:optTestSpacesImpl}*))
// [...] Definition of the object bilinearForm (see Section(*~\ref{sec:bilinform}*))
// Creation of the object systemAssembler
auto systemAssembler
   = make_DPG_SystemAssembler(nearOptTestSpaces, solutionSpaces, bilinearForm);
// Assemble the system Ax=f
MatrixType stiffnessMatrix;
VectorType rhsVector;
//[...] Define rhsFunction as a suitable lambda expression
systemAssembler.assembleSystem(stiffnessMatrix, rhsVector, rhsFunction);
// [...] Compute ``inverse(matrix)*rhsVector''
\end{c++}

\cpp{SystemAssembler} is also responsible for applying boundary conditions to the system. So far, only Dirichlet boundary conditions are implemented. To this end, first the degrees of freedom affected by the boundary condition are marked. Then, the boundary values are set to the corresponding nodes with the method \cpp{applyDirichletBoundarySolution}. For instance, if we are considering the transport problem of Section~\ref{sec:introTransport}, we need to set the degrees of freedom of $\theta$, that are in $\Gamma_-$, to $0$. For this, we mark the relevant nodes with the method \cpp{getInflowBoundaryMask} and store the information in a vector \cpp{dirichletNodesInflow}. Then we call \cpp{applyDirichletBoundarySolution} as we outline in the following listing. Note that the trial space associated to $\theta$ is required. Since, in our ordering, $\theta$ is our second unknown, we get its associated trial space with the command \cpp{std::get<1>(solutionSpaces)} (since \cpp{std::tuple} starts counting from 0).
\begin{c++}
  std::vector<bool> dirichletNodesInflow;
  BoundaryTools boundaryTools = BoundaryTools();
  boundaryTools.getInflowBoundaryMask(std::get<1>(solutionSpaces),
                                      dirichletNodesInflow,
                                      beta);    //mark affected degrees of freedom
  systemAssembler.applyDirichletBoundarySolution<1>
      (stiffnessMatrix,
       rhs,
       dirichletNodesInflow,
       0.);
\end{c++}

Finally, in certain types of problems, some degrees of freedom might be ill-posed. For example, in the transport case, the degrees of freedom corresponding to trial functions on faces aligned  with the flow direction will be weighted with 0 coefficients in the matrix. To address this issue, \cpp{SystemAssembler} provides several methods, of which the simplest is \cpp{defineCharacteristicFaces}.

\subsubsection{BilinearForm and InnerProduct}
\label{sec:bilinform}
As it follows from \eqref{eq:bilinSumK}, the bilinear form $b(\cdot, \cdot)$ can be decomposed into local bilinear forms $b_K(\cdot, \cdot)$. The \cpp{BilinearForm} class describes $b_K$ and provides access to the corresponding local matrices $A_K$ defined in \eqref{eq:AK} which are then used by the  \cpp{SystemAssembler} to assemble the global matrix $A$.

In our case, we view a bilinear form $b_K$ as a sum of what we will call elementary integral terms. By this we mean integrals over $K$ (or $\partial K$) which are a product of a test search function $v \in \bV_h$ (or its derivatives) and a trial function $u \in \bU_H$ (or its derivatives). Additionally, the product might also involve some given coefficient $c(x)$. For instance, in our transport equation (cf.~\eqref{eq:ultraweakdiscrete}),
\begin{equation}
\label{eq:bilinformexample}
b_K(u, v) = b_K((\varphi,\theta), v) \coloneqq
\underbrace{\int_K c v \varphi}_{Int_0} \underbrace{- \int_K \beta\cdot\nabla v \varphi}_{Int_1}
\underbrace{+ \sum_{K_h \in \Omega_h, K_h \subset K}\int_{\partial K_h} v \theta \beta \cdot n}_{Int_2},
\end{equation}
where we have omitted the tilde to ease notation here. 
Therefore the matrix $A_K = \sum_{i \in I} A_K^i$ can be computed as a sum of the matrices $A_K^i$ corresponding to the different elementary integrals $Int_i$, $i \in I$. Any of the elementary integrals can be expressed via the class \cpp{IntegralTerm} that we describe in Section~\ref{sec:integralterm}.

To create an object \cpp{bilinearForm} of the class \cpp{BilinearForm}, we call \cpp{make_BilinearForm} as follows.
\begin{c++}
auto bilinearForm = make_BilinearForm (testSearchSpaces, solutionSpaces, terms);
\end{c++}
The variables \cpp{testSearchSpaces} and \cpp{solutionSpaces} are the ones introduced in Section~\ref{sec:systemassembler} to represent $\bV_h$ and $\bU_H$. The object \cpp{terms} is a tuple of objects of the class \cpp{IntegralTerm}. Once that the object \cpp{bilinearForm} exists, a call to the method \cpp{getLocalMatrix} computes $A_K$ by iterating over all elementary integral terms and summing up their contributions $A_K^i$.

Let us now briefly discuss the class \cpp{InnerProduct}. Its aim is to allow the computation of the inner products associated to the Hilbert spaces $\bU$ and $\bV$. For this, we take advantage of the fact that  an inner product can be seen as a symmetric bilinear form $b(u,v)$ where $u$ and $v$ are both functions from some space. Hence, we can reuse the structure of \cpp{BilinearForm} for summing over elementary integral terms to define the class \cpp{InnerProduct}.
The construction of an \cpp{InnerProduct} is thus done with
\begin{c++}
auto innerProduct = make_InnerProduct (testSpaces, terms);
\end{c++}

\subsubsection{IntegralTerm}
\label{sec:integralterm}
An \cpp{IntegralTerm} represents an elementary integral over the interior of a cell $K$, over its faces $\partial K$ or even over faces of a partition of $K$. It expresses a product between a term related to a test function $v$ and a term related to a trial function $u$. Examples are $Int_0$, $Int_1$ and $Int_2$ from \eqref{eq:bilinformexample}.

The \cpp{IntegralTerm} is parametrized by two \cpp{size_t} that give the indices of the test and trial spaces that we want to integrate over.
Additionally we specify the type of evaluations used in the integral with a template parameter of type
\begin{c++}
enum class IntegrationType {
    valueValue,
    gradValue,
    valueGrad,
    gradGrad,
    normalVector,
    normalSign
};
\end{c++}
and the domain of integration with a template parameter of type
\begin{c++}
enum class DomainOfIntegration {
    interior,
    face
};
\end{c++}
If \cpp{integrationType} is of type
\cpp{IntegrationType::valueValue} or
\cpp{IntegrationType::normalSign}, the function \cpp{make_IntegralTerm} has to be called as follows:
\begin{c++}
auto integralTerm
   = make_IntegralTerm<lhsSpaceIndex, rhsSpaceIndex,
                       integrationType, domainOfIntegration>(c);
\end{c++}
where \cpp{c} is a scalar coefficient in front of the test space product and is of arithmetic type, e.\,g.\ \cpp{double}. The template parameter \cpp{domainOfIntegration} is one of the types from \cpp{DomainOfIntegration} and the parameters \cpp{lhsSpaceIndex} and \cpp{rhsSpaceIndex} refer, in this particular order, to the indices of test and trial space in their respective tuples of test and trial spaces. Note that the objects of the class \cpp{IntegralTerm} are not given the spaces themselves but only some indices referring to them. This is because the spaces are managed by the class \cpp{BilinearForm} (or \cpp{InnerProduct}) owning the \cpp{IntegralTerm}.

For other \cpp{integrationType}s, we also need to specify the flow direction \cpp{beta} by calling
\begin{c++}
auto integralTerm
   = make_IntegralTerm<lhsSpaceIndex, rhsSpaceIndex,
                       integrationType, domainOfIntegration>(c, beta);
\end{c++}
where \cpp{c} is again of arithmetic type and \cpp{beta} is of vector type, e.\,g.\ \cpp{FieldVector<double, dim>}. There has been some rudimentary work to support functions mapping coordinates to scalars or vectors for \cpp{c} and \cpp{beta}.

The \cpp{IntegralTerm} $Int_1$ from example~\eqref{eq:bilinformexample} can be created with
\begin{c++}
auto integralTerm
   = make_IntegralTerm<0, 0, IntegrationType::gradValue,
                       DomainOfIntegration::interior>(-1., beta);
\end{c++}
where the two zeroes are, in this particular order, the indices of test and trial space in their respective tuples of test and trial spaces.

The class \cpp{IntegralTerm} provides a method \cpp{getLocalMatrix} that computes its contribution $A_K^i$ to the local matrix $A_K$ and that is called by the \cpp{getLocalMatrix} method of \cpp{BilinearForm} or \cpp{InnerProduct}. To prevent runtime switches over the \cpp{IntegrationType} and \cpp{DomainOfIntegration}, we made them template parameters of \cpp{IntegralTerm} and use Boost Fusion\footnote{Fusion is a meta programming library and part of the C++ library collection Boost:~\url{http://www.boost.org/}} to easily handle compile time abstractions.

\subsection{Computing the optimal test space}
\label{sec:optTestSpacesImpl}
As described in Section~\ref{sec:localizationTheory}, for a given basis function $u_H^i$ of the trial space $\bU_H$, the corresponding near-optimal test function $\bar v^i=\sum_{K\in \Omega_H} \bar v^i_K \chi_K$ can be computed cell-wise. Indeed, one can find $\bar v^i_K$ by solving \eqref{eq:optBasisFiniteLocal} for every $K \in \Omega_H$. As a consequence, we can decompose the computation of $\bar v^i$ into the following steps (we will omit the index $i$ and call these functions $u_H,\ \bar v$ and $\bar v_K$ in the rest of this section):

\begin{itemize}
\item For a given cell $K\in \Omega_H$, let $\{z^j\}_{j=1}^M$ be a basis of $\bV_{h,K}$, the test-search space on cell $K$. In this basis, we can express $\bar v_K = \sum_{j=1}^M c_K^j z^j$ and find the vector of coefficients $c_K = ( c_K^j )_{j=1}^M$ as follows. Let us denote $B_K$ the $\bR^{M\times M}$ matrix with entries $(B_K)_{j,l}=\langle z^j,z^l \rangle_{\bV_K}$ for $1\leq j,l\leq M$ and $g_K \in \bR^{M}$ the vector with entries $g_K^j=b_{K} ( u_H,z^j),\ 1\leq j \leq M$. Then $c_K$ is the solution of the system
\begin{equation}
\label{eq:optBasisFiniteLocalMatrix}
B_K c_K = g_K.
\end{equation}
This task is done by the class \cpp{TestspaceCoefficientMatrix} for all basis functions $u_H$ with $\supp(u_H)\cap K \neq \emptyset$ (see Section~\ref{sec:TestspaceCoefficientMatrix} for more details).
\item Let $K_{ref}$ be the reference cell and $\{z^j_{ref}\}_{j=1}^M$ the basis functions of $\bV_{h,K}$ on $K_{ref}$. Having computed $c_K$, we define the corresponding local basis function of the near-optimal test space $\bar{v}_{K,ref} := \sum_{j=1}^M c_K^j z^j_{ref}$. This task is done by the class \cpp{OptimalTestLocalFiniteElement} (see Section~\ref{sec:OptimalTestLocalFiniteElement}).
\item Having computed $\bar{v}_{K,ref}$ and using the degree of freedom-handling from the global basis of the trial space $\bU_H$, we build  $\bar v_{ref}=\sum_{K\in \Omega_H} \bar v_{K,ref}$. This task is done by the class \cpp{OptimalTestBasis} (see Section~\ref{sec:OptimalTestBasis}). Note that the remaining mapping from $\bar v_{ref}$ to $\bar v$ has to be performed for any global basis. It is done with the help of geometry-information while assembling the system matrix $A$ and right hand side $F$.
\end{itemize}

\subsubsection{TestspaceCoefficientMatrix}
\label{sec:TestspaceCoefficientMatrix}
The computation of the coefficients $c_K$ is performed by the class \cpp{TestspaceCoefficientMatrix}
which has the bilinear form $b(\cdot, \cdot)$ and the inner product $\langle \cdot, \cdot \rangle_\mathbb{V}$ as template parameters. It has a method \cpp{bind(const Entity& e)}
in which it sets up and solves equation \eqref{eq:optBasisFiniteLocalMatrix} for all local basis functions $u_K$ of the trial space. Since the matrix $B_K$ is symmetric positive definite, the solution is determined via the Cholesky algorithm. The computed coefficients $c_K$ are saved in a matrix which can be accessed by the method \cpp{coefficientMatrix()}.

Furthermore, \cpp{TestspaceCoefficientMatrix} offers the possibility to save and reuse already computed coefficients. In case of constant parameters, equation \eqref{eq:optBasisFiniteLocalMatrix} only depends on the geometry of the cell $K$, so in case of a uniform grid, the coefficients are the same for all cells and thus need not be recomputed for every cell. To this end, in the current version, constant parameters are assumed and the geometry of the last cell is saved and compared to the geometry of the current cell. If they coincide, computation is skipped and the old coefficients are used. This is also helpful in cases with more than one or vector-valued test variables (see Section~\ref{sec:OptimalTestBasis}).

\subsubsection{OptimalTestLocalFiniteElement}
\label{sec:OptimalTestLocalFiniteElement}
The class \cpp{OptimalTestLocalFiniteElement} provides a local basis $\{ \bar{v}^i_{K,ref} \}_{i=1}^{\mathcal{N}}$ consisting of linear combinations $\bar{v}^i_{K,ref} = \sum_{j=1}^M (c_K^i)^j z^j_{ref}$ of a given local basis $\{z^j_{ref}\}_{j=1}^M$. Its constructor is called as follows
\lstset{escapeinside={(*}{*)},style=cppstyle}
\begin{c++}
OptimalTestLocalFiniteElement< D, R, d, TestSearchSpace>
                    (coefficientMatrix, //matrix containing the coefficients (*\smash{$c_K^i$}*)
                     testSearchSpace,   //given local basis (*\smash{$\{z_{ref}^j\}$}*)
                     offset=0)          //optional, default-value: 0    
\end{c++}
The template parameters are the type \cpp{D} used for domain coordinates, the type \cpp{R} used for function values, an integer \cpp{d} specifying the dimension of the reference element and the type of the given local basis \cpp{TestSearchSpace}. The optional parameter \cpp{offset} is used if the \cpp{coefficientMatrix} stores coefficients for multiple local bases at the same time. This happens for example for vector-valued test spaces. In this case, only the rows \cpp{offset} to \cpp{offset}$+$ size of \cpp{testSearchSpace} of the \cpp{coefficientMatrix} are taken into account.

\subsubsection{OptimalTestBasis}
\label{sec:OptimalTestBasis}
The near-optimal test space itself is implemented as a class called \cpp{OptimalTestBasis<TestspaceCoefficientMatrix, testIndex>}, which is compliant with the requirements for a global basis in the dune-module \dunefunctions~\citep{EGMS2015}. In the current implementation, \cpp{OptimalTestBasis} describes only one scalar variable. If there are more test variables, several \cpp{OptimalTestBasis} are needed. In this case, \cpp{testIndex} specifies the index of the component of the near-optimal test space which is described.

Like all global bases in \dunefunctions , the \cpp{OptimalTestBasis}\ provides a \cpp{LocalIndexSet}\ for the mapping of local to global degrees of freedom as well as a \cpp{LocalView}\ to provide access to all local basis functions whose support has non-trivial intersection with a given element. Since every optimal test basis function corresponds to a basis function in the trial spaces, the mapping of local to global degrees follows directly from the corresponding \cpp{LocalIndexSets}\ of the trial spaces. If there is more than one trial space, the global degrees of freedom are ordered by trial space, that is the first global degrees of freedom are those for the first trial space, the next ones are for the second trial space and so on. The \cpp{LocalView}\ of the \cpp{OptimalTestBasis}\ provides access to the correct \cpp{OptimalTestLocalFiniteElement}. To this end, in the method \cpp{bind(const Element& e)}, the coefficients for the local near-optimal test basis are computed by binding the \cpp{coefficientMatrix} to the element \cpp{e} and the corresponding \cpp{OptimalTestLocalFiniteElement} is constructed.

Note that for test spaces with more than one component, the computation of the coefficients for the local near-optimal test basis in general cannot be seperated for the different components. That is why the \cpp{TestspaceCoefficientMatrix} is implented as a seperate class so that several near-optimal test spaces can share one \cpp{TestspaceCoefficientMatrix} and the computed coefficients can be reused if the \cpp{TestspaceCoefficientMatrix} is bound several times to the same element by different near-optimal test bases.

The following lines of code show how to create an \cpp{OptimalTestBasis} consisting of two variables for a given \cpp{BilinearForm} and \cpp{InnerProduct}.
\begin{c++}
  TestspaceCoefficientMatrix testspaceCoefficientMatrix(bilinearForm, innerProduct);
  Functions::OptimalTestBasis<TestspaceCoefficientMatrix, 0>  feBasisTest0(testspaceCoefficientMatrix);
  Functions::OptimalTestBasis<TestspaceCoefficientMatrix, 1>  feBasisTest1(testspaceCoefficientMatrix);
\end{c++}

\subsection{A posteriori error estimators}
\label{sec:apost}
To compute the residual
\begin{equation*}
\Vert f - \cB u_H \Vert_{\bV'}
=
\sup_{v\in \bV} \frac{\Vert f(v) - b(u_H,v) \Vert_{\bV}}{\Vert v \Vert_{\bV}},
\end{equation*}
we exploit once again the product structure of $\bV$ and use the fact that
\begin{equation*}
\Vert f - \cB u_H \Vert^2_{\bV'}
= \sum_{K\in \Omega_H} \Vert r_{K}(u_H,f) \Vert^2_{\bV_K'}
= \sum_{K\in \Omega_H} \Vert R_{K}(u_H,f) \Vert^2_{\bV_K}
\end{equation*}
where $r_{K}$ is the cell-wise residual. $R_{K}$ is the Riesz-lift of $r_{K}$ in $\bV_K$ so it is the solution of
\begin{equation}
\label{eq:resGalerkin}
\left< R_{K}(u_H,f),v \right>_{\bV_K}
=
b(u_H,v)-f(v), \quad \forall v \in \bV_K.
\end{equation}
Since \eqref{eq:resGalerkin} is an infinite-dimensional problem, we project the Riesz-lift $R_K$ to a finite-dimensional subspace $\overline\bV_K$ of $\bV_K$, obtaining an approximation $\overline R_K$.
This in turn gives the a posteriori error estimator
\begin{equation}
\label{eq:apost}
\Vert \overline R(u_H,f) \Vert_{\bV}
\coloneqq
\left(
\sum_{K\in \Omega_H} \Vert \overline R_K(u_H,f) \Vert^2_{\bV_K}
\right)^{1/2}.
\end{equation}
An appropriate choice of the a posteriori search space $\overline\bV_K$ depends on the problem and is crucial to make $\Vert \overline R_K \Vert_{\bV_K}$ be good error indicators.

In \dunedpg, the computation of the a posteriori estimator \eqref{eq:apost} is handled by the class \cpp{ErrorTools}. The following lines of code compute \eqref{eq:apost} for the solution \cpp{u_H} of a problem with bilinear form \cpp{bilinearForm}, inner product  \cpp{innerProduct} and right hand side \cpp{rhsVector}. 
\begin{c++}
ErrorTools errorTools = ErrorTools();
double aposterioriErr = errorTools.aPosterioriError(bilinearForm,innerProduct,u_H,rhsVector); 
\end{c++}
The object \cpp{bilinearForm} is of the type \cpp{BilinearForm} described above. It has to be created with an object \cpp{testSpace} associated to the a posteriori search space $\overline \bV_H$. The same applies for \cpp{innerProduct}, which is of type \cpp{InnerProduct}.

\section{Numerical Example: Implementation of Pure Transport in \dunedpg}
\label{sec:numEx}
\begin{figure}[h]
\centering
\begin{subfigure}[b]{0.49\textwidth}
\includegraphics[width=\textwidth]{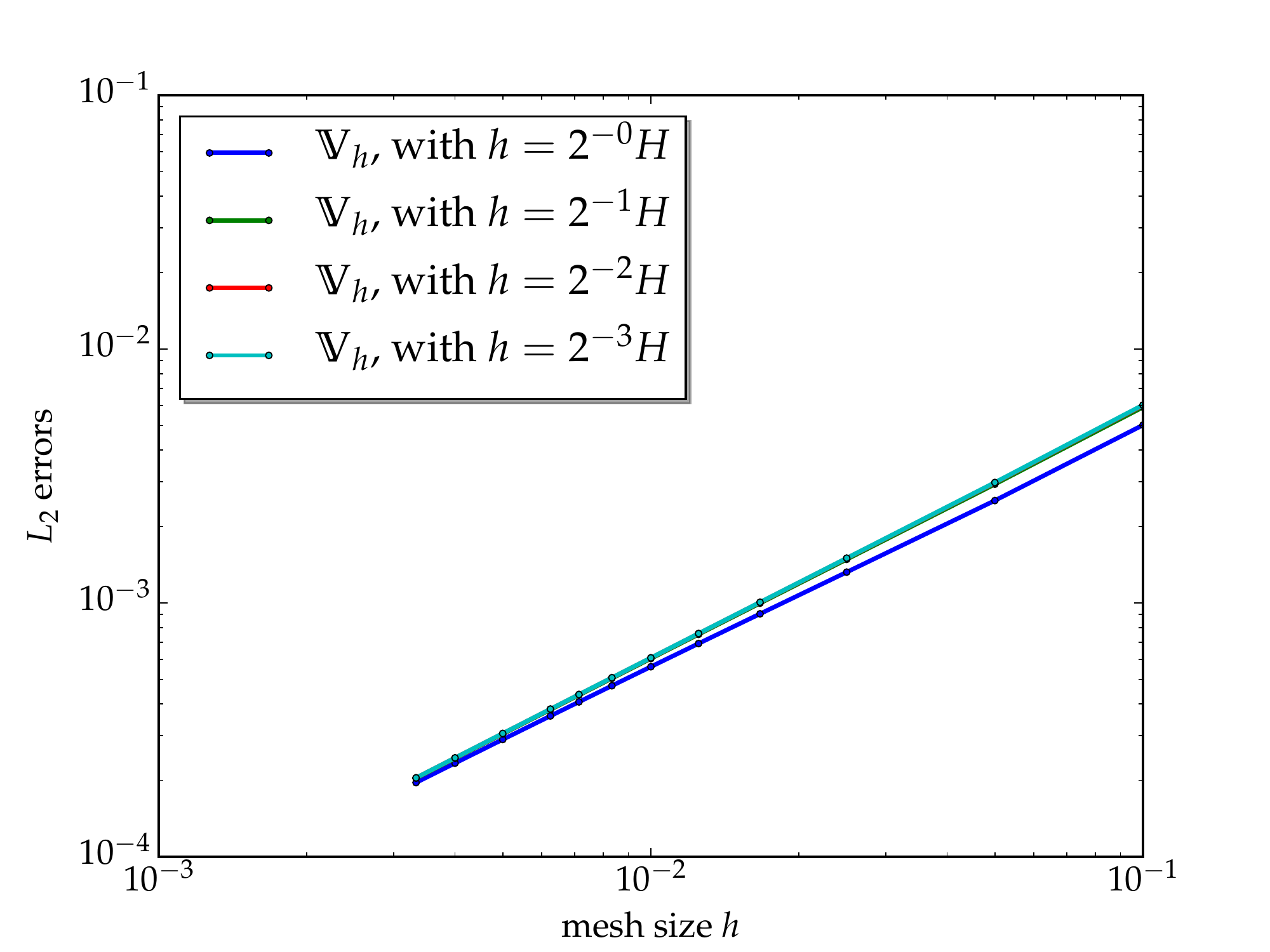}
\caption{$L_2$ error of $\varphi$}
\label{fig:error}
\end{subfigure}
~
\begin{subfigure}[b]{0.49\textwidth}
\includegraphics[width=\textwidth]{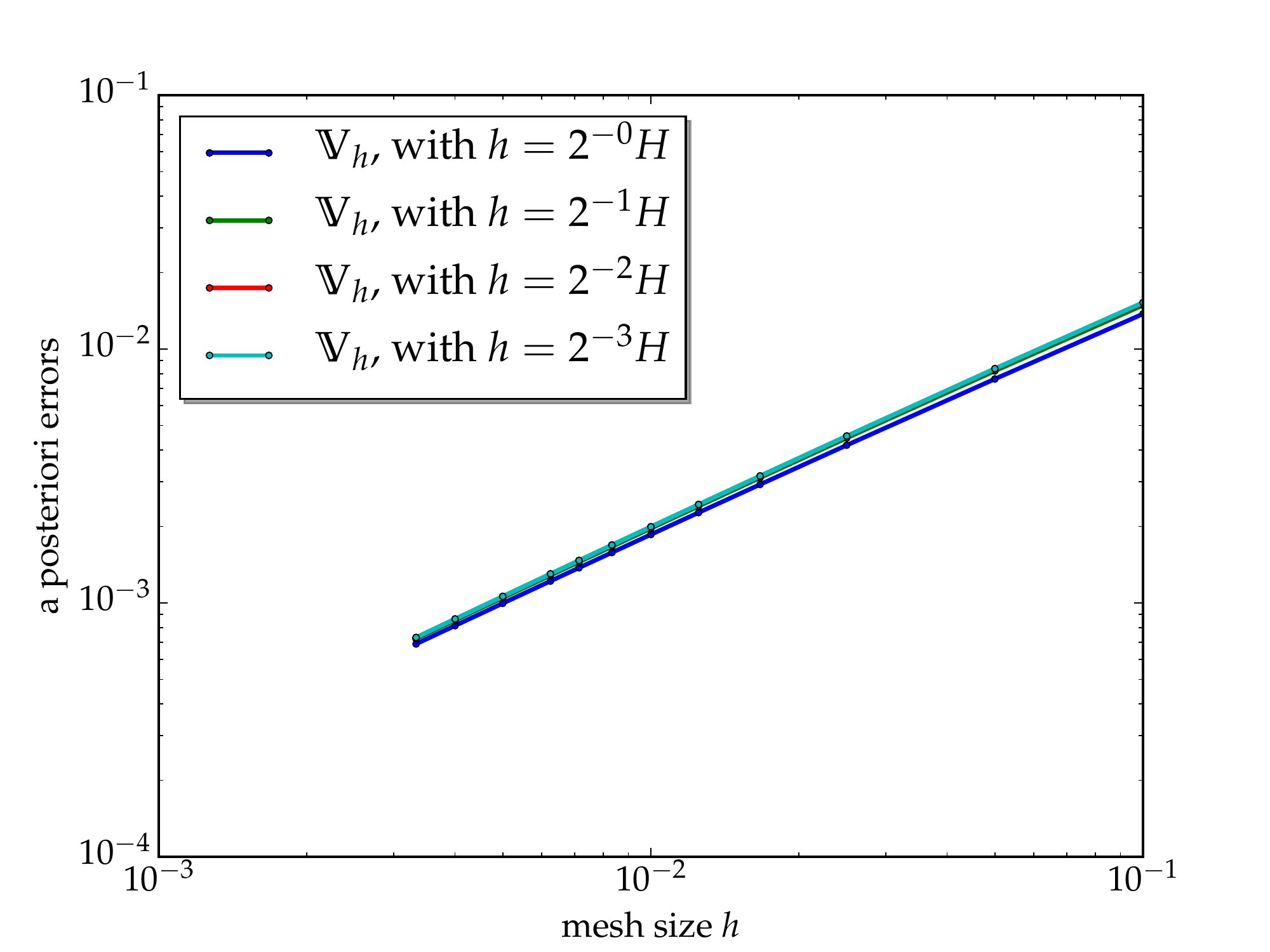}
\caption{a posteriori error of $u=(\varphi, \theta)$}
\label{fig:error_aposteriori}
\end{subfigure}
\\
\begin{subfigure}[b]{0.49\textwidth}
\includegraphics[width=\textwidth]{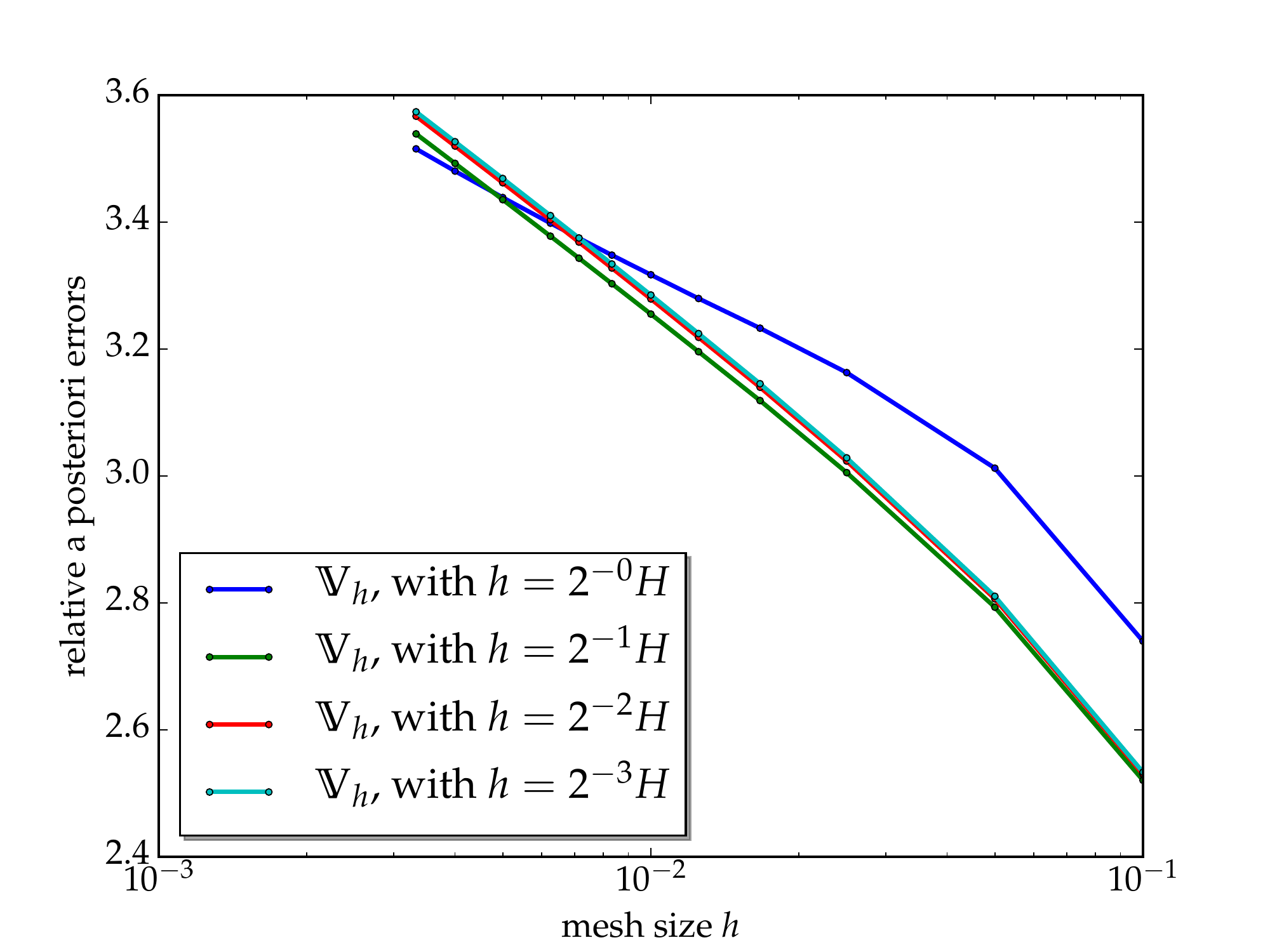}
\caption{relative a posteriori error
         $\frac{e_{\mathrm{aposteriori}}}{e_{\mathrm{exact}}}$}
\label{fig:error_apost_rel}
\end{subfigure}
~
\begin{subfigure}[b]{0.49\textwidth}
 \includegraphics[width=\textwidth, clip=true, trim=15cm 0 0 2cm]{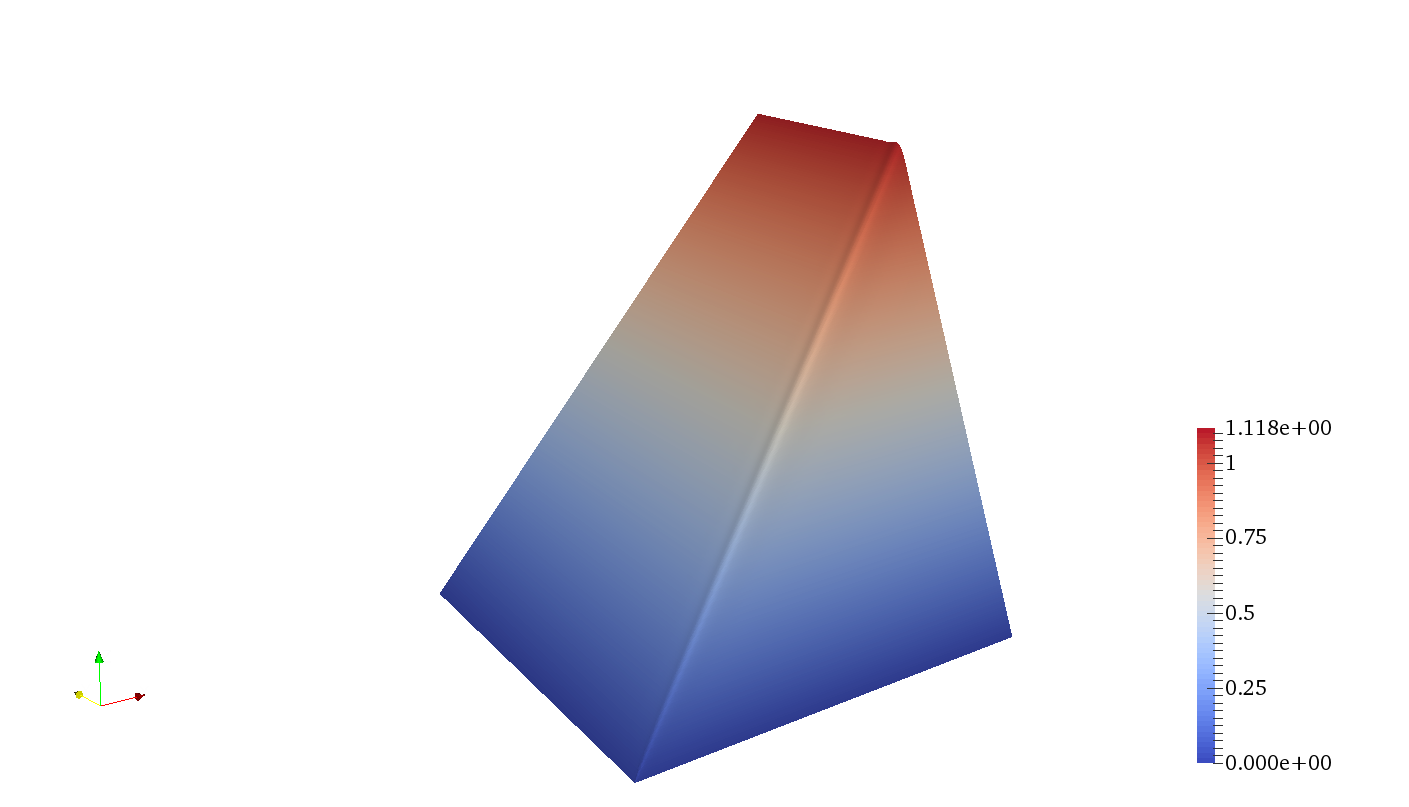}
\caption{Solution $\varphi$.}
\label{fig:exactSolution}
\end{subfigure}
\caption{$L_2$ error and a posteriori error estimator of numerical solutions}
\end{figure}

As a simple numerical example, we solve the transport problem \eqref{eq:pureTransport} with 
\begin{align*}
c&=0, \\
\beta&=(\cos(\pi/8),\, \sin(\pi/8)), \\
f&=1.
\end{align*}
As Figure~\ref{fig:exactSolution} shows, the exact solution $\varphi$ describes a linear ramp starting at $0$ in each point of the inflow boundary $\Gamma_-$ and increasing with slope 1 along the flow direction $\beta$. There is a kink in the solution starting in the lower left corner of $\Omega$ and propagating along $\beta$.

For the numerical solution, we let $\Omega_H$ be a partition of $\Omega$ into uniformly shape regular triangles. $\Omega_h$ is a refinement of $\Omega_H$ to some level $\ell \in \bN_0$ such that $h=2^{-\ell}H$.

With $\bU_H$ and $\bV_h$ defined as in \eqref{Uh} and \eqref{testsearch} with $m=2$, we compute $u_H = (\varphi_H,\theta_H) \in \bU_H$ by solving the ultra-weak variational formulation \eqref{eq:finiteVFtrans}. We investigate convergence in $H$ of the error $\Vert \varphi - \varphi_H \Vert_{L_2(\Omega)}$. We also evaluate the a posteriori estimator $\Vert \overline R(u_H,f) \Vert_{\bV}$ when the components $\overline R_K(u_H,f)$ are computed with a subspace $\overline \bV_K$ of polynomials of degree 5, $\forall K \in \Omega_H$.

Regarding the error $\Vert \varphi - \varphi_H \Vert_{L_2(\Omega)}$, as Figure~\ref{fig:error} shows, we observe linear convergence as $H$ decreases. This is to be expected since the polynomial degree to compute $\varphi_H$ is 1. The figure also shows that the refinement level $\ell$ of the test-search space $\bV_h$ has essentially no impact on the behavior of the error.

Regarding the behavior of the a posteriori estimator $\Vert \overline R(u_H,f) \Vert_{\bV}$, it is possible to see in Figures~\ref{fig:error_aposteriori} and \ref{fig:error_apost_rel} that the quality of $\Vert \overline R(u_H,f) \Vert_{\bV}$ slightly degrades as $H$ decreases in the sense that, as $H$ decreases, $\Vert \overline R(u_H,f) \Vert_{\bV}$ represents the error $\Vert \varphi - \varphi_H \Vert_{L_2(\Omega)}$ less and less faithfully. An element that might be playing a role is that $\Vert \overline R_{K}(u_H,f)\Vert_{\bV}$ is not exactly an estimation of $\Vert \varphi - \varphi_H \Vert_{L_2(\Omega)}$, but of the error including also $\theta_H$, namely $\Vert u - u_H\Vert_{\bU} = \Vert (\varphi,\theta) - (\varphi_H,\theta_H)\Vert_{\bU}$.

\section{Conclusion And Future Work}
In Section~\ref{sec:foundationsDPG}, we gave a short overview of the DPG method.
We then introduced our \dunedpg~library in Section~\ref{sec:dunedpg}, documenting the internal structure and showing how to use it to solve a given PDE.
Finally, we showed some numerical convergence results computed for a problem with well-known solution. This allowed us to compare our a posteriori estimators to the real $L_2$ error of our numerical solution. 
As a next step we want to implement adaptive mesh refinements that would be driven by our local a posteriori error indicators.

Finally, we want to improve our handling of vector valued problems with one notable example being first order formulations of convection--diffusion problems.
With our current \cpp{std::tuple} of spaces structure used throughout the code, we have to implement vector valued spaces by adding the same scalar valued space several times.
With the \dunetypetree\ library from \citet{Muething2015} we can handle vector valued spaces much more easily, as has already been shown in \dunefunctions.
This will result in mayor changes in our code, but will probably allow us to replace our dependency on Boost Fusion with more modern C++11 constructs.
In the long run, we hope that this would give us increased maintainability and decreased compile times in addition to the improvements in the usability of vector valued problems.
This is aligned with our long-term goal of making \dunedpg\ a flexible building block for constructing DPG solvers for a large range of different problem types.

\section*{Acknowledgments}
We thank O.~Sander for his introduction to the \dune~library and his guidance in understanding it. We also thank W.~Dahmen for introducing us to the topic of DPG. Finally, O.~Mula is indebted to the AICES institute of RWTH~Aachen for hosting her as a postdoc during 2014--2015 which is the period in which large parts of \dunedpg\ were developed.

\bibliographystyle{abbrvnat}
\bibliography{literature}
\end{document}